\documentclass[11pt,fleqn]{amsart}
\usepackage{amsmath}
\usepackage{amssymb}
\usepackage[all]{xy}           
\usepackage{amssymb}           
\usepackage{hyperref}
\usepackage{eucal}
\usepackage{graphicx}
\usepackage{epsfig}
\usepackage[usenames,dvipsnames]{color}
\usepackage{charter}
\usepackage{color}
\numberwithin{equation}{section}

\newtheorem{definition}{Definition}[section]

\newtheorem{theorem}[definition]{Theorem}
\newtheorem{proposition}[definition]{Proposition}

\newtheorem{remarkth}[definition]{Remark}

\newenvironment{remark}{\begin{remarkth}\upshape}{\hfill$\diamond$\end{remarkth}}

\begin{document}

\title{Poisson structures on smooth 4--manifolds}

\author{Luis C. Garc\'ia-Naranjo}

\address{Departamento de Matem\'aticas y Mec\'anica \\ IIMAS-UNAM \\ Apdo Postal 20-726, Mexico City, 01000, Mexico}

\email{luis@mym.iimas.unam.mx}

\author{Pablo Su\'arez-Serrato}

\address{Instituto de Matem\'aticas - Universidad Nacional Aut\'onoma de M\'exico\\Circuito Exterior, Ciudad Universitaria\\Coyoac\'an, 04510\\Mexico City\\Mexico}

\email{pablo@im.unam.mx}

\author{Ram{\'o}n Vera}

\address{The Pennsylvania State University \\Mathematics Department
 \\ University Park, State College,  \\ PA 16802, USA}

\email{rvera.math@gmail.com}

\begin{abstract}
We show that every closed oriented smooth 4-manifold admits a complete singular Poisson structure in each homotopy class of maps to the 2-sphere. The  rank  of this structure is 2 outside a small singularity set, which consists of finitely many circles and isolated points. The Poisson bivector vanishes on the singularities, where we give its local form explicitly.
\end{abstract}


\maketitle

\section{Introduction and results}

What is the most prevalent geometric structure available on a smooth 4-manifold?  Symplectic 4-manifolds have been given a precise description in terms of Lefschetz pencils, beginning with the seminal work of S.K. Donaldson in \cite{D99}. The notion of symplectic Lefschetz fibrations was then extended to include singularities along circles in \cite{ADK05}. These fibrations are now known as {\em broken Lefschetz fibrations}. The existence of such structures on closed\footnote{By {\it closed} manifold we mean a compact manifold without boundary.} smooth oriented 4-manifolds has been shown in \cite{AC08, B08, L09}. Every smooth, oriented, closed 4-manifold admits a broken Lefschetz fibration. Our contribution here is to indicate the existence of a Poisson structure of rank 2 whose symplectic leaves (as fibres) and singularities coincide precisely with those of a broken Lefschetz fibration. Moreover, to explicitly present the local forms of the Poisson bivector and the associated symplectic forms on the neighbourhoods of the singularities. We thus obtain:

\begin{theorem}\label{main-theo}
Let $X$ be a closed oriented smooth 4-manifold. Then on each homotopy class of map from $X$ to the 2-sphere there exists a complete Poisson structure of rank 2 on $X$ whose associated Poisson bivector vanishes only on a finite collection of circles and isolated points.  
\end{theorem}

A proof of this Theorem, given in section \ref{proof-main}, follows from a known construction in the field of Poisson 
geometry, described in detail in section 2.1. (We thank P. Frejlich for indicating this construction to us). In the first version of this paper we provided an alternative proof by giving local descriptions of the bivector and a consistent glueing.

%

This construction equips the fibres (both singular and regular) of the broken Lefschetz fibration 
with a symplectic form away from the singularities. These symplectic forms
approach infinity close to the singularities (since the Poisson structure approaches zero). 
The precise rate at which this happens is explained in section \ref{S:Symp-Struct-Sing}.
This contrasts with near-symplectic constructions on broken Lefschetz fibrations, where the opposite phenomenon occurs. That is, this closed 2-form, called near-symplectic, is non-degenerate outside a set of singular circles, where it vanishes.

Notice that every oriented 2-manifold admits a symplectic structure and hence a Poisson structure. For dimension three, it was shown in \cite{IM03} that every closed oriented 3-manifold admits a regular rank 2 Poisson structure. In dimension 4 only symplectic manifolds can have constant rank 4 Poisson structures. As we pointed out, there is a close relationship between symplectic 4-manifolds and Lefschetz fibrations. So our construction provides families of Poisson bi-vectors of rank 2 on general closed oriented 4-manifolds that in some sense possess the smallest singular sets possible.

For example, on the standard four-sphere $S^4$ the only constant rank Poisson structure must have rank 0. This follows since $S^4$ does not admit a symplectic structure, so it can not admit a Poisson structure of constant rank 4. As the Euler class of its tangent bundle is zero, it does not split as the direct sum of two orientable 2-plane bundles, which would be the case if it were to admit a constant rank 2 Poisson structure. In contrast, through our result we can now endow $S^4$ with a Poisson structure of generic rank 2 for every homotopy class of maps to $S^2$.

Another noteworthy consequence of Theorem \ref{main-theo}, due to M. Kontsevich \cite{K03}, is that every closed oriented smooth 4-manifold admits a deformation quantization on each homotopy class of map to the 2-sphere. 


Our paper aims to bring the fields of Poisson geometry and smooth 4-manifolds closer together. So we provide thorough introductions to Poisson manifolds in section \ref{ss:Poisson} and to broken Lefschetz fibrations in section \ref{ss:BLF}. Finally in section \ref{ss:examples} we present various examples. 

{\bf Acknowledgements:} We thank Jaume Amor\'os, Paula Balseiro, Henrique Bursztyn,  
Pedro Frejlich, David Iglesias, Yank{\i}  Lekili,  Juan Carlos Marrero, Eva Miranda, Edith Padr\'on, Gerardo Sosa and Jonathan Williams for interesting conversations and correspondence. PSS thanks CONACyT Mexico and PAPIIT UNAM for supporting various research activities. RV thanks Jos{\'e} Seade for all his support during his stay at the Instituto de Matem{\'a}ticas, UNAM. We are also grateful to an anonymous referee for numerous valuable remarks.

\section{Definitions}

\subsection{Poisson manifolds}\label{ss:Poisson}
We give a short review of the standard definitions and we present some results from the field of Poisson geometry that will be needed in the sequel. The reader can find more details in the books \cite{V94, DZ05}.
\begin{definition}
\label{D:Poisson}
 A {\em Poisson bracket} (or a {\em Poisson structure}) on a smooth manifold $M$ is a bilinear operation
 $\{\cdot , \cdot \}$ on the set $C^\infty(M)$ of real valued smooth functions on $M$ such that 
 \begin{enumerate}
\item[(i)] $( C^\infty(M) , \{\cdot , \cdot \})$ is a Lie algebra.
\item[(ii)] $\{\cdot , \cdot \}$ is a derivation in each factor, that is,
\begin{equation*}
\{gh, k\}=g\{h, k\}+ h\{g, k\}
\end{equation*}
for any $g,h,k\in  C^\infty(M)$.
\end{enumerate}
\end{definition}
A manifold $M$ endowed with a Poisson bracket is called a {\em Poisson manifold}. The class of Poisson morphisms on $M$
is given by the so-called {\em Poisson maps} defined below.

\begin{definition}
\label{D:PoissonMap}
A smooth map $\varphi:M\to M$ on a Poisson manifold $M$ is called a {\em Poisson map} if
\begin{equation*}
\{g\circ \varphi , h \circ \varphi \} = \{ g, h \} \circ \varphi
\end{equation*}
for all $g, h \in C^\infty (M)$.
\end{definition}

The most basic and fundamental example of a non-trivial Poisson manifold is a symplectic manifold $(M,\omega)$.
The bracket on $M$ is defined by
\begin{equation*}
\{g,h\}=\omega(X_g,X_h).
\end{equation*}
Recall that the Hamiltonian vector field $X_h$ is defined through the relation ${\bf i}_{X_h}\omega =dh$ and similarly
for $X_g$. The Jacobi identity for the bracket follows from the property of $\omega$ being closed.

The derivation property (ii) in Definition \ref{D:Poisson} allows one to extend the notion of Hamiltonian vector fields 
beyond the symplectic setting.  Given a function $h\in C^\infty(M)$  we associate to it
the {\em Hamiltonian vector field} $X_h$, that is defined as the following derivation on $C^\infty(M)$
\begin{equation*}
X_h(\cdot )=\{\cdot , h \}.
\end{equation*}
It is a simple exercise to show that in a symplectic manifold the Poisson and symplectic definitions of Hamiltonian vector fields
are consistent.

Another consequence of (ii) is that the bracket $\{g,h\}$ only depends on the first derivatives of $g$ and $h$.
So we may think of the bracket as defining a contravariant antisymmetric 2-tensor $\pi$ on $M$ such that 
\begin{equation}\label{eq:bracket-bivector}
\{g,h\}=\pi(dg,dh).
\end{equation}
Note that $\pi$ is a   section  of $\Lambda^2 TM$, i.e. $\pi$ is a bivector field. It is common to refer to $\pi$  as the {\em Poisson tensor}. In local coordinates $(x^1, \dots , x^n)$ we can represent	
\begin{equation*}
\pi(x)=\frac{1}{2}\sum_{i,j=1}^n\pi^{ij}(x)\frac{\partial}{\partial x^i}\wedge \frac{\partial}{\partial x^j},
\end{equation*}
where $\pi^{ij}(x)=\{x^i,x^j\}=-\{x^j,x^i\}$. 

The Jacobi identity for the bracket implies that $\pi$ satisfies an integrability  condition which in local coordinates
is a system of first order semilinear partial differential equations for $\pi^{ij}(x)$. It can also be expressed intrinsically 
as $[\pi,\pi]=0$, where $[\cdot,\cdot]$ is the Schouten-Nijenhuis bracket of multivector fields \cite{V94}. We mention that
this bracket is a local operator so the Jacobi identity is a local condition on $\pi$.

In the sequel we will often say that the Poisson bivector $\pi$ is the Poisson structure on the manifold $M$.

To a general bivector $\pi$ in $M$ (that is not necessarily Poisson) 
one can associate a bundle map $\mathcal{B}:T^*M\to TM$ defined by its action on covectors by the rule $$\mathcal{B}_q(\alpha_q)(\cdot)=\pi_q(\cdot ,\alpha_q)$$ where
$\alpha_q\in T_q^*M$ and $q\in M$. If $\pi$ is Poisson, then $X_h=\mathcal{B}(dh)$.
 
 The {\em rank} of a bivector $\pi$  at a point $q\in M$ is defined to be the rank of $\mathcal{B}_q:T^*_qM\to T_qM$. 
In local coordinates it is the rank of the matrix $\pi^{ij}(x)$. The image of $\mathcal{B}_q$ is a 
subspace $D_q\subset T_qM$, and the collection of these subspaces as $q$ varies on $M$ defines the so-called
{\em characteristic distribution} of $\pi$.  Note that, at every point, the  rank of the 
characteristic distribution is even and
coincides with the  rank of $\pi$, and that the characteristic distribution is smooth but may not be regular. 
If $\pi$ is Poisson,  the rank of $\pi$ at $q\in M$ is called the {\em  rank of the Poisson structure at $q$}.

The celebrated  {\em Symplectic Stratification Theorem} states that the characteristic distribution
of any Poisson tensor $\pi$ is integrable\footnote{Here and in what follows, we use the term {\em integrable} in the sense of integrability of singular distributions as is made clear in the Stefan-Sussman Theorem. See e.g. \cite{DZ05}. } . Denote by  $\Sigma_q$  the even dimensional leaf of the
corresponding foliation of $M$ passing through the point $q\in M$. One can characterize $\Sigma_q$
as the set of points in $M$ that can be joined
with $q$ with a piecewise smooth curve, each of which is a trajectory of a locally defined Hamiltonian vector
field. The  Symplectic Stratification Theorem also guarantees that $\Sigma_q$ is an even dimensional immersed submanifold 
of $M$ that carries a symplectic structure $\omega_{\Sigma_q}$. Moreover, the theorem asserts that the Poisson bracket 
$\{ \cdot , \cdot \}_{\omega_{\Sigma_q}}$ on $\Sigma_q$ induced by $\omega_{\Sigma_q}$ coincides with the ``restriction" of
$ \{\cdot , \cdot \}$ to $\Sigma_q$. More precisely, if $g,h\in C^\infty(\Sigma_q)$ then, for any $z\in \Sigma_q$ we have
\begin{equation*}
\{g , h \}_{\omega_{\Sigma_q}}(z)=\{\tilde g , \tilde h \}(z)
\end{equation*}
where $\tilde g , \tilde h$ are arbitrary smooth extensions of $g, h$ to $M$.

Note that $T_q\Sigma_q=D_q$. Thus, if $u_q, v_q \in T_q\Sigma_q$
there exist co-vectors $\alpha_q, \beta_q\in T^*_qM$ that correspondingly map to $u_q$ and $v_q$ by $\mathcal{B}_q$.
One has 
\begin{equation}
\label{E:Symp-form-gen}
\omega_{\Sigma_q}(q)(u_q, v_q)=\pi_q(\alpha_q, \beta_q)=\langle \alpha_q,  v_q \rangle=-\langle \beta_q,  u_q \rangle.
\end{equation}

The set $\Sigma_q$ is called {\em the symplectic leaf through $q$} and  in this way we obtain 
a {\em foliation of $M$ by symplectic leaves}.

In general, the rank of a Poisson structure is not constant and hence
the leaves of the symplectic foliation will have different dimension. It is common to say that the {\em rank of a Poisson structure} is the maximal rank of the structure at each point in $M$, and to refer to the points of lower rank as {\em singular points}. A Poisson structure with no singular points is called {\em regular}.


A very special situation occurs if the rank of the characteristic distribution of a bivector  is less than or equal to two.
\begin{proposition}
\label{P:rank2}
If $\pi$ is a bivector field on $M$ whose characteristic distribution is integrable and has rank less than or equal to two
at each point, then $\pi$ is Poisson.
\end{proposition}
\begin{proof}
One only needs to verify that the Jacobi identity $[\pi, \pi]=0$ holds. 

At a point $q$ where the rank of $\pi$  is zero, one has $\pi_q=0$, and from the local expression for $[\pi, \pi]$ (see e.g. \cite{DZ05}) it immediately follows that $[\pi, \pi](q)=0$.

On the other hand, since the rank of $\pi$ is at most 2, the set of points where the rank of $\pi$ is 2 is open in $M$.
Hence, if the rank of $\pi$ at $q\in M$ is 2,  there exists a neighbourhood $U\subset M$  that contains $q$ 
such that $\pi|_U$ has constant rank 2. 

By integrability of the characteristic distribution of $\pi$,  the bivector $\pi|_U$ defines a foliation of $U$ by 
2-dimensional leaves.  Each leaf $\Sigma$ is equipped with a 2-form $\omega_\Sigma$ defined as follows. By definition of the
characteristic distribution, for $q\in \Sigma$ and $v_q, w_q\in T_p\Sigma$ there exist covectors $\alpha_q, \beta_q\in T_q^*U$ with
\begin{equation*}
v_q=\mathcal{B}_q(\alpha_q), \qquad w_q=\mathcal{B}_q(\beta_q).
\end{equation*}
Then $\omega_\Sigma(q)(v_q, w_q):=\pi_q(\alpha_q, \beta_q)$. These 2-forms are nondegenerate and vary smoothly from leaf to leaf. 
Moreover, since the dimension of the leaves is  two, they are also closed. The leafwise 2-forms being closed is equivalent to the Jacobi identity (see Chapter 2 
of \cite{V94} for more details). Therefore, $[\pi|_U, \pi|_U]=0$. But this implies that $[\pi, \pi]|_U=0$ since the Jacobi identity
is a local condition.
\end{proof}
As a consequence of the above proposition we also have
\begin{proposition}
\label{P:Rank2conformal}
Let $\pi$ be a Poisson structure on $M$ whose rank at each point is less than or equal to two. Then $\pi_1:=k\pi$ is also a Poisson
structure where $k\in C^\infty(M)$ is an arbitrary non-vanishing function.
\end{proposition}
\begin{proof}
Since $k$ does not vanish, the characteristic distribution of $\pi_1$ coincides with that of $\pi$, and is thus integrable. The
result now follows from Proposition \ref{P:rank2}.
\end{proof}

\begin{definition}
\label{D:Casimir}
Let $M$ be a Poisson manifold. A function $h\in C^\infty(M)$ is called a {\em Casimir} if $\{h,g\}=0$ for 
every $g\in C^\infty(M)$.
\end{definition}

Suppose that $M$ is an orientable $n$-manifold and that let $\mu$ be an orientation of $M$ (a non-vanishing $n$-form on $M$).
Let $F_1, \dots, F_{n-2}\in C^\infty(M)$ and consider the bivector 
on $M$ defined by the relation
\begin{equation*}
\pi(\alpha, \beta) \mu:=k\,\alpha\wedge \beta\wedge dF_1\wedge\dots \wedge dF_{n-2}
\end{equation*}
for 1-forms $\alpha$, $\beta$ on $M$, and 
where $k\in C^\infty(M)$ is a non-vanishing function. 

We can give a complete description of the characteristic distribution of $\pi$.
Let $F:=(F_1, \dots ,F_{n-2}):M\to  \mathbb{R}^{n-2}$. If $q\in M$ is a critical point of $F$, meaning that $dF_1(q), \dots , dF_{n-2}(q)$
are linearly dependent, then $D_q=\{0\}$. Otherwise, $D_q$ is tangent to the level set of $F$ containing $q$ which (at least locally)
is a 2-dimensional submanifold of $M$. Hence, the characteristic distribution of $\pi$ is integrable and its leaves are:
\begin{enumerate}
\item The 2-dimensional  leaves $F^{-1}(y)$ where $y\in  \mathbb{R}^{n-2}$ is a regular value of $F$,
\item the 2-dimensional  leaves  $F^{-1}(y)\setminus C$  where $y\in  \mathbb{R}^{n-2}$ is a critical value of $F$ and 
$C:=\{\mbox{Critical Points of $F$}\}\subset M$,
\item the zero dimensional leaves consisting of each point in $C$.
\end{enumerate}

By Proposition \ref{P:rank2} we conclude that $\pi$ is Poisson. The specific form of the bracket is determined from the
following expression:
\begin{equation}
\label{E:DefintermsofCasimirs}
\{g,h\} \mu:=k\,dg\wedge dh\wedge dF_1\wedge\dots \wedge dF_{n-2}.
\end{equation}
It is straightforward to see that $F_1, \dots, F_{n-2}$ are Casimirs.

Formula \eqref{E:DefintermsofCasimirs} appeared for the first time in \cite{Dam89} (attributed to H. Flaschka and T. Ratiu) and admits the
following generalization 
that will be key for our purposes.
\begin{theorem}
\label{T:Const-Poisson}
Let $M$ be an orientable $n$-manifold, $N$ an orientable $n-2$ manifold, and $f:M\to N$ a smooth map.
Let $\mu$ and $\Omega$ be orientations of $M$ and $N$ respectively. The bracket on $M$ defined by
\begin{equation}
\label{E:Def-Intrinsic}
\{g,h\}\mu=k\,dg\wedge dh \wedge f^*\Omega
\end{equation}
where $k$ is any non-vanishing function on $M$ is Poisson. Moreover, its symplectic leaves are
\begin{enumerate}
\item[(i)] the 2-dimensional leaves $f^{-1}(s)$ where $s\in N$ is a regular value of $f$,
\item[(ii)] the 2-dimensional leaves $f^{-1}(s)\setminus C$ where  $s\in N$ is a singular value of $f$ and 
$C:=\{\mbox{Critical Points of $f$}\}\subset M$,
\item[(iii)] the 0-dimensional leaves corresponding to each element in $C$.
\end{enumerate}
\end{theorem}
\begin{proof}
In local coordinates \eqref{E:Def-Intrinsic} takes the form \eqref{E:DefintermsofCasimirs} where $(F_1,\dots, F_{n-2})$
are the components of the coordinate representation of $f$ (the function  $k$ need not be the same in both expressions but
remains non-vanishing). The proof follows from the discussion given before equation \eqref{E:DefintermsofCasimirs}.
\end{proof}

\begin{remark}
\label{R:Arbit-k}
Note that a different choice of the forms $\mu$ and $\Omega$ in \eqref{E:Def-Intrinsic} results in the appearance of a positive factor
that can be absorbed into the function $k$. Hence, the construction outlined in the  theorem gives a family of Poisson
 brackets on $M$, parametrized by non-vanishing functions $k\in C^\infty(M)$, 
 whose symplectic foliation satisfies items (i), (ii), (iii). This freedom in the choice
 of $k$ should be expected in view of Proposition \ref{P:Rank2conformal}.
\end{remark}

We mention that other generalizations of \eqref{E:DefintermsofCasimirs} appear in \cite{DamP12} to construct
Poisson brackets of rank greater than 2 with prescribed Casimir functions.

%

\begin{definition}
A Poisson manifold $M$ is said to be complete if every Hamiltonian
vector field on $M$ is complete.  
\end{definition}

Notice that $M$ is complete if and only if every symplectic leaf is bounded in the
sense that its closure is compact. 

\subsection{Broken Lefschetz fibrations}\label{ss:BLF}
Before defining the concept of broken Lefschetz fibrations we will say a word about Lefschetz fibrations.  A {\it Lefschetz fibration} on a simply connected 4-manifold $X$ is a smooth map $f\colon X\rightarrow S^2$ whose generic fibre is a surface.  The map $f$ is allowed to have isolated critical points, known as Lefschetz singularities, which are modeled in local complex coordinates by $f\colon (z_1, z_2) \rightarrow z_{1}^{2} + z_{2}^{2}$.  Regular fibres are smooth, but singular fibres present an isolated nodal singularity.  

A Lefschetz pencil on a 4-manifold $X$ is a map $f\colon  X \setminus B \rightarrow  S^2$,  
which is not 
defined at a finite number of base points  $\{b_1, \dots , b_m \} =  B$.  Around each base point, $f$ is modeled in 
local complex coordinates by $f\colon (z_1, z_2) \mapsto z_1 / z_2$.  Alternatively, thinking of $S^2$ as $\mathbb{CP}^1$,
then $f\colon (z_1, z_2) \mapsto [z_1 : z_2]$.  The fibres of $f$ are punctured surfaces, to which one adds the base points 
to obtain closed surfaces, called the fibres of the pencil.  Near a base point, a piece of a fibre looks like the slicing of 
$\mathbb{C}^2$ into complex planes passing through the origin.  If one blows up a Lefschetz pencil at all its base points, then one obtains a Lefschetz fibration. 
\begin{figure}[h]
    \centering
    \includegraphics[width=0.8\textwidth]{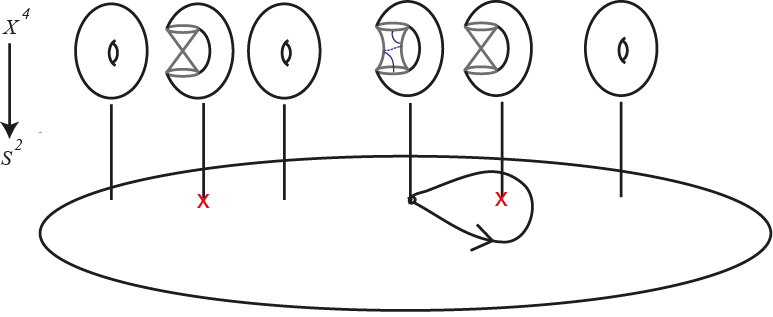}
    \caption{Example of a Lefschetz fibration}
    \label{fig:LF}
\end{figure}

In Figure \ref{fig:LF}, we can see an example of a Lefschetz fibration with generic fibre the 2-torus. As the regular fibre approaches the singular Lefschetz point, the vanishing cycle shrinks to a point.   Over the critical values the fibre presents an isolated singularity.  However, outside the singularities all fibres are $T^{2}$.  The figure also depicts the monodromy phenomena.  Start with a generic fibre over a regular value.  Go around a singular point in a closed loop. This action accounts for a positive Dehn twist along the vanishing cycle on the regular fibre (see the lectures by Seidel and Smith in \cite{ACMSSST}).

Lefschetz fibrations are closely related to symplectic structures. The work of Donaldson and Gompf shows that Lefschetz fibrations are in direct correspondence with symplectic 4-manifolds.  Symplectic 4-manifolds admit Lefschetz fibrations, up to blow ups.  Furthermore, given a Lefschetz fibration on a smooth oriented 4-manifold with a suitable cohomology class, it is posible to construct a symplectic structure on the total space with symplectic fibres \cite{D99}. 

Auroux, Donaldson and Katzarkov introduced a second type of singularity to the Lefshetz fibrations to study a larger class of smooth 4-manifolds that are non-symplectic. These mappings are known as {\it broken Lefschetz fibrations} or BLF.  This terminology was introduced by Perutz in \cite{P07-2}. By a BLF, we understand a submersion $f\colon X \rightarrow S^2$ with two types of singularities: isolated points (Lefschetz singularities) and circles (indefinite folds).

\begin{definition}
\label{D:BLF}
On a smooth, closed 4-manifold $X$, a {\it broken Lefschetz fibration} or {\it BLF} is a smooth map $f\colon X \rightarrow S^2$ that is a submersion outside the singularity set.  Moreover, the allowed singularities are of the following type:
\begin{enumerate}
\item {\it Lefschetz} singularities:  finitely many points  $$C_{L}= \{ p_1, \dots , p_k\} \subset X,$$  which are locally modeled by complex charts
$$ \mathbb{C}^{2} \rightarrow \mathbb{C}  ,  \quad \quad (z_1, z_2) \mapsto z_{1}^{2} + z_{2}^{2},$$
\item {\it indefinite fold} singularities, also called {\it broken}, contained in the smooth embedded 1-dimensional submanifold $\Gamma \subset X \setminus  C_{L}$, which are locally modelled by the real charts
$$ \mathbb{R}^{4} \rightarrow \mathbb{R}^{2} ,  \quad \quad  (t,x_1,x_2,x_3) \mapsto (t, - x_{1}^{2} + x_{2}^{2} + x_{3}^{2}).$$
\end{enumerate}
We will denote $C=C_{L} \cup \Gamma$ the set of singular points of $f$.
\end{definition}
%


The term indefinite in (ii) refers to the fact that the quadratic form $- x_{1}^{2} + x_{2}^{2} + x_{3}^{2}$ is neither negative nor positive definite. In the language of singularity theory, these subsets are known as fold singularities of corank 1. Since $X$ is closed, $\Gamma$ is homeomorphic to a collection of disjoint circles.  
 For this reason, throughout this work, we will often refer to $\Gamma$ as {\it singular circles}.  On the other hand, we can
 only assert that $f(\Gamma)$ is  a union of immersed curves. In particular, the images of the components of $\Gamma$ need not be disjoint, and the image of each component can self-intersect. 

On a normal neighbourhood of a circle within $\Gamma$  we  have a splitting of $\mathbb{R}^{3}\times S^1 \to S^1$ into a real line bundle over $S^1$ and a rank 2  bundle over $S^1$. There are two homotopy classes of such splittings, orientable and non-orientable ones \cite{Ho204}. Thus, there are two models, a product $S^1\times B^3$ which is orientable, and a non-orientable model given by the quotient of $S^1\times B^3$ by an involution that reverses the orientations on both summands of the splitting, for example $(\theta,x_1,x_2,x_3) \mapsto (\theta+\pi, -x_1, -x_2, x_3)$. 

\begin{figure}
    \centering
    \includegraphics[width=0.6\textwidth]{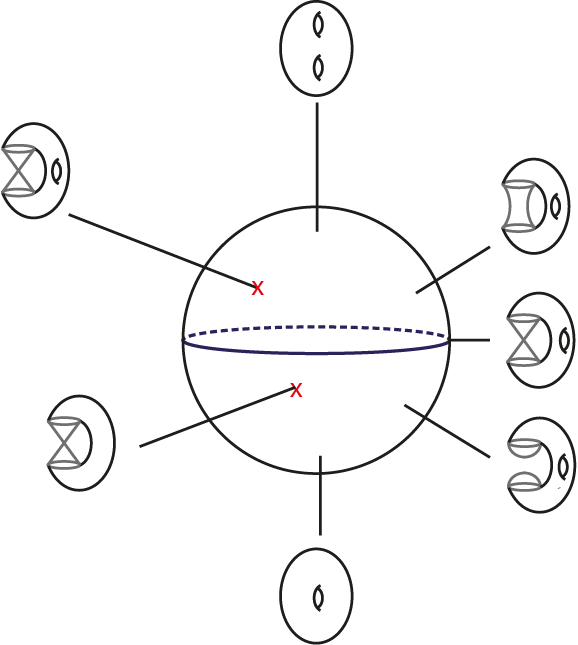}
    \caption{A diagram depicting a broken Lefschetz fibration with one circle of folds and two Lefschetz singularities}
    \label{fig:BLF}
\end{figure}

Figure \ref{fig:BLF} depicts an example of a broken Lefschetz fibration.  This example considers only one singular circle. The image of this circle under $f$ is shown at the equator of the 2-sphere.  Over the northern hemisphere of $S^{2}$ the fibres are genus 2 surfaces. Crossing the image of the singular circle amounts to a change in the topology of the fibre.  The fibres on the southern hemisphere are tori.  On each hemisphere there is one Lefschetz singular point, where the fibre has an isolated nodal singularity. 

\begin{remark} In general fold singularities can generate disconnected regular fibres. However, it was shown in \cite{B08} that, after a homotopy, the fibres of a BLF may be assumed to be connected. \end{remark}

\begin{remark}A priori, the singularities of a BLF can appear in a complicated way. The circles of folds could intersect each other or the Lefschetz points could lie between the circles. Nevertheless, it is possible to obtain a BLF with a simple representation. These are called {\it simplified broken Lefschetz fibrations}, which is a BLF $f\colon X\rightarrow S^{2}$ with only one circle of indefinite folds whose image is on the equator, and with all critical Lefschetz values lying on one hemisphere. Corollary 1 of \cite{W10} provides an existence result for a simplified broken Lefschetz fibration in any homotopy class of maps from a smooth $4$--manifold to the 2-sphere, which in turn relates to results from  \cite{L09} and \cite{B09}.
\end{remark}

%
%
%
%
%
%
\subsection{Proof of  Theorem \ref{main-theo}}\label{proof-main}
\begin{proof}
Let $X$ be a closed smooth oriented and connected 4-manifold, and $f\colon X \to S^2$ a broken Lefschetz fibration. 
The hypothesis of Theorem  \eqref{T:Const-Poisson} are met with $X=M$ and $N=S^2$. So, for any
non-vanishing $k\in C^\infty(X)$, the expression \eqref{E:Def-Intrinsic} defines a
Poisson structure $\pi$ on $X$ with the required properties. The structure is complete since $X$ is compact. 
\end{proof}


\begin{remark}
\label{R:Arbit-k-X}
As in Remark \ref{R:Arbit-k}, we stress that the freedom in the choice of $k$ is a natural consequence of the
rank of $\pi$ being less than or equal to 2 (Proposition \ref{P:Rank2conformal}).
\end{remark}

\begin{remark}
We point out that the above proof implies that the {\em wrinkled fibrations} of \cite{L09} can also be given a compatible Poisson structure. The associated local normal forms for these Poisson bivectors and their deformations, as well as the growth rates of the symplectic forms on the fibres, will be given in a forthcoming note.
\end{remark}

\section{Local expressions for the bracket and the induced symplectic forms on the fibres near the singularities.}\label{S:Symp-Struct-Sing}

We will now construct explicit expressions for the Poisson structure and the corresponding symplectic forms in the
vicinity of the singularities of a   BLF $f\colon X \to S^2$.
All of the expressions that we give depend on the particular choice of the non-vanishing function $k\in C^\infty(X)$ (see
Remarks  \ref{R:Arbit-k} and  \ref{R:Arbit-k-X}).

\subsection{Local expressions around Lefschetz singularities.}
\label{SS:IsolatedSingularity}

As  described
in item (i) of Definition \ref{D:BLF}, the coordinate representation of $f$ near a Lefschetz singularity $p\in X$ is given by $F:=(F_1,F_2):B^4\to B^2$ where
\begin{equation*}
\begin{split}
F_1(x_1,y_1,x_2,y_2)&=x_1^2-y_1^2+x_2^2-y_2^2, \\ F_2(x_1,y_1,x_2,y_2)&=2(x_1y_1+x_2y_2),
\end{split}
\end{equation*}
and $B^4\subset \mathbb{R}^4$ and $B^2\subset \mathbb{R}^2$ are neighbourhoods around $0$. Here $F_1$ and $F_2$ are the real and imaginary parts of the parametrization function $\varphi:\mathbb{C}^2\to \mathbb{C}$ given by
\begin{equation*}
\varphi(z_1,z_2)=z_1^2+z_2^2,
\end{equation*}
where $z_j=x_j+iy_j$.

Applying formula \eqref{E:DefintermsofCasimirs} yields the local expression for $\pi$ near a Lefschetz singularity.
%
%
%
%
%
\begin{equation}
\label{E:pi_C}
\begin{split}
\pi=&k(x_1,y_1,x_2,y_2) \left ( (x_2^2+y_2^2)\frac{\partial}{\partial x_1}\wedge \frac{\partial}{\partial y_1} +(x_1^2+y_1^2)\frac{\partial}{\partial x_2}\wedge \frac{\partial}{\partial y_2}  \right .\\ & + (-y_1y_2-x_1x_2) \frac{\partial}{\partial x_1}\wedge \frac{\partial}{\partial y_2}
+(x_1x_2+y_1y_2)\frac{\partial}{\partial y_1}\wedge \frac{\partial}{\partial x_2}\\ & \left . +(-x_1y_2+y_1x_2)\frac{\partial}{\partial y_1}\wedge \frac{\partial}{\partial y_2}+(-x_1y_2+y_1x_2)\frac{\partial}{\partial x_1}\wedge \frac{\partial}{\partial x_2}\right ).
\end{split}
\end{equation}
As usual, $k$ is any smooth non-vanishing function.


Now we give an expression for the symplectic form on the fibres near $p$.  Consider a non-zero point $q=(x_1,y_1,x_2,y_2)\in B^4$ and assume for the moment that $x_1^2+y_1^2\neq 0$.
 The symplectic leaf $\Sigma_q$ through $q$  is a  level set of $F$ (taking away $0\in B^4$ if necessary).
A simple calculation shows that the vectors
\begin{equation}
\label{E:Vecuv-piC}
\begin{split}
u_q&=\frac{1}{\sqrt{(x_1^2+y_1^2)(x_1^2+y_1^2+x_2^2+y_2^2)}}\left(-(x_1x_2+y_1y_2)\frac{\partial}{\partial x_1}\right . \\
& \qquad \qquad \qquad \left .- (x_1y_2-x_2y_1)\frac{\partial}{\partial y_1}+(x_1^2+y_1^2)\frac{\partial}{\partial x_2}\right ),\\
v_q&=\frac{1}{\sqrt{(x_1^2+y_1^2)(x_1^2+y_1^2+x_2^2+y_2^2)}}\left((x_1y_2-x_2y_1)\frac{\partial}{\partial x_1}\right . \\
& \qquad \qquad \qquad \left .-(x_1x_2+y_1y_2) \frac{\partial}{\partial y_1}+(x_1^2+y_1^2)\frac{\partial}{\partial y_2}\right ),
\end{split}
\end{equation}
are annihilated by $dF_1(q)$ and $dF_2(q)$. Hence
they are tangent to $\Sigma_q$ at $q$. Moreover, they are orthonormal with respect to the euclidean metric
\begin{equation*}
ds^2=dx_1^2+dy_1^2+dx_2^2+dy_2^2
\end{equation*}
on $B^4$. Using \eqref{E:pi_C}
one can check that $\mathcal{B}_p(\alpha_p)=u_p$, where
$$\alpha_p=\frac{dy_2}{ k(x_1,y_1,x_2,y_2)\sqrt{(x_1^2+y_1^2)(x_1^2+y_1^2+x_2^2+y_2^2)}}.$$
Therefore, in view of \eqref{E:Symp-form-gen} we have
\begin{equation}
\label{E:Symp-form-pi-C}
\omega_{ \Sigma_p}(p)(u_p,v_p)=\langle \alpha_p, v_p \rangle = \frac{1}{k(x_1,y_1,x_2,y_2)(x_1^2+y_1^2+x_2^2+y_2^2)}.
\end{equation}
We can now prove the following.
\begin{proposition}
\label{P:symp-structure-piC}
Let $q\in B^4\setminus \{{ 0 }\}$. The symplectic form induced by $\pi$ on the symplectic leaf $\Sigma_q$ through 
$q=(x_1,y_1,x_2,y_2)$ at the
point $q$ is given by
\begin{equation*}
\omega_{\Sigma_q}(q)=\frac{1}{k(x_1,y_1,x_2,y_2)(x_1^2+y_1^2+x_2^2+y_2^2)}\,\omega_{Area}(q)
\end{equation*}
where $\omega_{Area}$ is the area form
on $\Sigma_q$ induced by the euclidean metric on $B^4$.
\end{proposition}
\begin{proof}
 Assume that $x_1^2+y_1^2\neq 0$. Given that $u_p, v_p$ 
defined by \eqref{E:Vecuv-piC} are orthonormal with respect to the euclidean metric we
have
\begin{equation*}
\omega_{Area}(p)(u_p,v_p)=1.
\end{equation*}
The result now follows from \eqref{E:Symp-form-pi-C}. Indeed, since the symplectic leaf $\Sigma_q$ is 2-dimensional, the 2-forms $\omega_{\Sigma_q}$ and $\omega_{Area}$
evaluated at $q$ must be proportional.

If $x_1^2+y_1^2= 0$ but $x_2^2+y_2^2\neq 0$, the same argument can be applied to the vectors
\begin{equation*}
u_q= \frac{\partial}{\partial x_{1}}, \qquad v_q =\frac{\partial}{\partial y_{1}}.
\end{equation*}

In this case $\mathcal{B}_q(\alpha_q)=u_q$, where
$$\alpha_q=\frac{dy_1}{k(0,0,x_2,y_2)( x_2^2+y_2^2)}.$$
Therefore,
\begin{equation*}
\omega_{\Sigma_q}(q)(u_q,v_q)=\langle \alpha_q, v_q\rangle =  \frac{1}{k(0,0,x_2,y_2)(x_2^2+y_2^2)}.
\end{equation*}
\end{proof}

%
%
%
%
%
%
%

\subsection{Local expressions near singular circles.}
\label{SS:SingularCircle}
Consider a singular circle within $\Gamma\subset X$.  In this case we can give an expression for the Poisson tensor $\pi$ not only
in a neighbourhood of a point in the circle, but on the whole normal bundle of the circle.

First we consider the case when the normal bundle of the singular circle is orientable.
As described after Definition  \ref{D:BLF}, the  model for this bundle is $S^1\times B^3$ with $f$ represented as
$F=(F_1,F_2):S^1\times B^3\to B^2$ given by
\begin{equation*}
\begin{split}
F_1(\theta, x_1,x_2,x_3)&=\theta, \\
F_2(\theta, x_1,x_2,x_3)&=-x_1^2+x_2^2+x_3^2.
\end{split}
\end{equation*}
A straightforward application of \eqref{E:DefintermsofCasimirs} yields the following expression for $\pi$  in the normal
bundle of a singular circle (valid if such bundle is orientable).
\begin{equation} \label{eq:pi-gamma}
\pi=k(\theta,x_1,x_2,x_3)\left (x_1 \frac{\partial}{\partial x_2} \wedge \frac{\partial}{\partial x_3} + x_2 \frac{\partial}{\partial x_1} \wedge \frac{\partial}{\partial x_3} - x_3 \frac{\partial}{\partial x_1} \wedge \frac{\partial}{\partial x_2} \right ),
\end{equation}
where $k$ is a non-vanishing function.
The above tensor can be interpreted as a multiple of a linear Poisson structure in $\mathbb{R}^3$. Hence, up to the factor $k$, it is dual to the Lie algebra structure
of real dimension 3 possessing the following commutation relations between the basis elements $e_1, e_2, e_3$:
\begin{equation*}
[e_1,e_2]=-e_3, \qquad [e_2,e_3]=e_1, \qquad [e_1,e_3]=e_2 .
\end{equation*}
This Lie algebra  is isomorphic to  $\frak{sl}(2,\mathbb{R})$. Therefore, in the vicinity of
a singular circle whose normal bundle is orientable, $\pi$  is proportional to the product Poisson structure of  $S^1$ equipped with the 
zero Poisson structure and $B^3$ equipped with the Lie-Poisson structure of $\frak{sl}(2,\mathbb{R})^*$ with
an appropriate basis identification.

When the normal bundle of the singular  circle is not orientable, we still have a coordinate description in terms of the quotient by the action of 
\begin{equation}
\label{E:defiota}
\iota\colon (\theta,x_1,x_2,x_3) \to (\theta+\pi, -x_1, -x_2, x_3).
\end{equation}
\begin{proposition}
The involution $\iota$ defined on $S^1\times B^3$ by \eqref{E:defiota}
 is a  Poisson map for the bracket  defined by \eqref{eq:pi-gamma} provided that $k\circ\iota=k$.
\end{proposition}
\begin{proof}
Let $g,  h\in C^\infty( S^1\times B^3)$. We have
\begin{align}
\lbrace g, h \rbrace (\theta,x_1,x_2,x_3)&=k(\theta,x_1,x_2,x_3)\left (x_1\left ( \frac{\partial g}{\partial x_2} \frac{\partial h}{\partial x_3} - \frac{\partial g}{\partial x_3} \frac{\partial h}{\partial x_2} \right ) \right .
\nonumber
\\
&+x_2\left ( \frac{\partial g}{\partial x_1} \frac{\partial h}{\partial x_3} - \frac{\partial g}{\partial x_3} \frac{\partial h}{\partial x_1} \right )
\nonumber
\\ 
&\left . - x_3\left ( \frac{\partial g}{\partial x_1} \frac{\partial h}{\partial x_2} - \frac{\partial g}{\partial x_2} \frac{\partial h}{\partial x_1} \right )
\right )
\nonumber
\end{align}
with all derivatives evaluated at $(\theta,x_1,x_2,x_3)$. Therefore,
\begin{align}
\{g, h\}\circ \iota (\theta,x_1,x_2,x_3)
=&k(\theta, x_1, x_2, x_3) \left ( -x_1\left ( \frac{\partial g}{\partial x_2} \frac{\partial h}{\partial x_3} - \frac{\partial g}{\partial x_3} \frac{\partial h}{\partial x_2} 
\right ) \right .
\nonumber
\\ 
&-x_2\left ( \frac{\partial g}{\partial x_1} \frac{\partial h}{\partial x_3} - \frac{\partial g}{\partial x_3} \frac{\partial h}{\partial x_1} \right )
\nonumber
\\
& \left . - x_3\left ( \frac{\partial g}{\partial x_1} \frac{\partial h}{\partial x_2} - \frac{\partial g}{\partial x_2} \frac{\partial h}{\partial x_1} \right )
\right )
\nonumber
\end{align}
with all derivatives evaluated at $\iota (\theta,x_1,x_2,x_3)=(\theta +\pi, -x_1,-x_2,x_3)$ and where we have used
$k\circ\iota=k$.

On the other hand,
\begin{align}
\{g\circ \iota , h\circ \iota \} (\theta,x_1,x_2,x_3)&=k(\theta, x_1, x_2, x_3) \left ( x_1\left ( \frac{\partial (g\circ \iota)}{\partial x_2} \frac{\partial (h\circ \iota)}{\partial x_3} - \frac{\partial  (g\circ \iota)}{\partial x_3} \frac{\partial  (h\circ \iota)}{\partial x_2} \right )  \right .
\nonumber
\\
 &+x_2\left ( \frac{\partial  (g\circ \iota)}{\partial x_1} \frac{\partial  (h\circ \iota)}{\partial x_3} - \frac{\partial  (g\circ \iota)}{\partial x_3} \frac{\partial  (h\circ \iota)}{\partial x_1} \right ) 
\nonumber
\\ 
&\left . - x_3\left ( \frac{\partial  (g\circ \iota)}{\partial x_1} \frac{\partial  (h\circ \iota)}{\partial x_2} - \frac{\partial  (g\circ \iota)}{\partial x_2} \frac{\partial  (h\circ \iota)}{\partial x_1} \right ) \right )
\nonumber
\end{align}
with all derivatives evaluated at $(\theta,x_1,x_2,x_3)$. However, by the chain rule
\begin{equation*}
\begin{split}
\frac{\partial (g\circ \iota)}{\partial x_i}(\theta,x_1,x_2,x_3)&=- \frac{\partial g}{\partial x_i}(\theta +\pi, -x_1,-x_2,x_3), \qquad i=1,2, \\
\frac{\partial (g\circ \iota)}{\partial x_3}(\theta,x_1,x_2,x_3)&= \frac{\partial g}{\partial x_3}(\theta +\pi, -x_1,-x_2,x_3),
\end{split}
\end{equation*}
and similarly for $h$. Hence $\{g, h\}\circ \iota=\{g\circ \iota , h\circ \iota \}$ as required.
\end{proof}

The above proposition implies  that if $k\circ \iota=k$, the expression  \eqref{eq:pi-gamma} 
 descends to the quotient of the product $S^1\times B^3$ by the action of $\iota$. Therefore, 
it also gives a valid expression for the bracket in a normal bundle of a singular circle if the latter is non-orientable.

In analogy with Proposition \ref{P:symp-structure-piC} we have
\begin{proposition}
\label{P:symp-structure-piGamma}
Let $q\in S^1\times (B^3\setminus \{{ 0}\})$. The symplectic form induced by $\pi$ on the symplectic leaf $\Sigma_q$ through 
$q=(\theta,x_1,x_2,x_3)$ at the
point $q$ is given by
\begin{equation*}
\omega_{\Sigma_q}(q)=\frac{1}{k(\theta,x_1,x_2,x_3)\sqrt{x_1^2+x_2^2+x_3^2}}\,\omega_{Area}(q)
\end{equation*}
where $\omega_{Area}$ is the area form
on $\Sigma_q$ induced by the metric $ds^2=d\theta^2+dx_1^2+dx_2^2+dx_3^2$ on $S^1\times B^3$.
\end{proposition}
\begin{proof}
We proceed analogously as we did in Section \ref{SS:IsolatedSingularity}. In this case the symplectic
leaf $\Sigma_q$ consists of points having the same $\theta$ coordinate as $q$ and belonging to the same 
level set of $F_2$ (removing the point $(\theta,0,0,0)$ if necessary).

First assume $x_2^2+x_3^2\neq 0$. The tangent vectors at $q$ 
\begin{equation*}
\begin{split}
u_q&=\frac{1}{\sqrt{x_2^2+x_3^2}}\left ( -x_3\frac{\partial}{\partial x_2} +x_2\frac{\partial}{\partial x_3}\right ), \\
v_q&=\frac{1}{\sqrt{x_2^2+x_3^2}\sqrt{x_1^2+x_2^2+x_3^2}}\left (  (x_2^2+x_3^2)\frac{\partial}{\partial x_1}  \right . \\ 
& \qquad \qquad  \qquad \qquad\left . +x_1x_2
\frac{\partial}{\partial x_2}+x_1x_3  \frac{\partial}{\partial x_3}\right ),
\end{split}
\end{equation*}
are annihilated by $dF_1(q)$ and $dF_2(q)$ so they are tangent to $\Sigma_q$. Moreover,
they are orthonormal with respect to the metric $ds^2=d\theta^2+dx_1^2+dx_2^2+dx_3^2$.
Hence
\begin{equation*}
\omega_{Area}(q)(u_q,v_q)=1.
\end{equation*}
On the other hand we have $v_q=\mathcal{B_q}(\beta_q)$ where
\begin{equation*}
\beta_q=\frac{x_3\, dx_2 -x_2\, dx_3}{k(\theta,x_1,x_2,x_3)\sqrt{x_2^2+x_3^2}\sqrt{x_1^2+x_2^2+x_3^2}}.
\end{equation*}
Therefore, in view of \eqref{E:Symp-form-gen} we have
\begin{equation*}
\omega_{ \Sigma_q}(q)(u_q,v_q)=-\langle \beta_q, u_q \rangle = \frac{1}{k(\theta,x_1,x_2,x_3)\sqrt{x_1^2+x_2^2+x_3^2}}.
\end{equation*}
The result now follows since $\omega_{ \Sigma_q}(q)$ and $\omega_{Area}(q)$ are proportional
as $\Sigma_q$ has dimension 2.

If $x_2=x_3=0$ and $x_1\neq 0$, the same argument can be applied with the vectors
\begin{equation*}
u_q=\frac{\partial}{\partial x_2}, \qquad v_q=\frac{\partial}{\partial x_3}.
\end{equation*}
This time $\beta_q=-\frac{dx_2}{x_1k(\theta,x_1,x_2,x_3)}$ so we get
\begin{equation*}
\omega_{ \Sigma_q}(q)(u_q,v_q)=\frac{1}{k(\theta,x_1,x_2,x_3)x_1}.
\end{equation*}
The sign ambiguity that arises when $x_1<0$ is fixed  by changing the orientation of  $\omega_{Area}(q)$.
\end{proof}

\section{Examples}\label{ss:examples}

\subsection{Near-symplectic 4-manifolds}\label{nssection}

Let $X = S^1\times Y^3$, where $Y$ is a closed Riemannian 3-manifold equipped with a circle-valued Morse function $f\colon Y \rightarrow S^{1}$ with indefinite type singularities. Locally, such functions have the same parametrization as real-valued Morse functions. That is, on a neighbourhood $U$ of a critical point of a Morse function  we have $f_{U} \colon (x_1, x_2, x_3) \mapsto (\pm x^{2}_{1} \pm x^{2}_{2} \pm x^{2}_{3})$.  A Morse function is called of {\it indefinite type} if it has no maximum nor minimum. Consider the Morse 1-form $df \in \Omega^{1}(Y)$ and denote by $t$  the angle parameter of $S^1$. 

The manifold $X$ constructed above is also an example of a near-symplectic manifold.  On a smooth, closed, oriented 4-manifold, a closed 2-form $\omega$  is {\it near-symplectic} if either it is non-degenerate, or it vanishes on a finite collection of disjoint circles.  This weaker condition makes near-symplectic 4-manifolds more abundant than symplectic ones.  If a smooth, oriented 4-manifold $X$ is compact and a certain topological condition holds ($b_{2}^{+}>0$), then there is a near-symplectic form $\omega$ on $X$ \cite{Ho04}.

\begin{definition}[\cite{ADK05}]\label{nsdef}
Let $X$ be a smooth oriented 4-manifold. Consider a closed 2-form $\omega \in \Omega^2(X)$ such that $\omega^2\geq 0$ and such that $\omega_p$ only has rank 4 or rank 0 at any point $p \in X$, but never rank 2. The form $\omega$ is called {\it near-symplectic}, if for every $p\in X$, either
\begin{enumerate}
 \item $\omega_p^2>0$, or
 \item $\omega_p = 0$, and $\textnormal{Rank}(\nabla \omega_p) = 3$, where $\nabla\omega_p\colon T_pX \rightarrow \Lambda^2 T^{*}_{p} X$ denotes the intrinsic gradient of $\omega$.
\end{enumerate}
\end{definition}

 It was shown in \cite{P07} that the zero set $Z_{\omega} = \lbrace p \in X \mid \omega_p = 0 \rbrace$ of a near-symplectic form $\omega \in \Omega^2 (X) $ is a smooth 1-dimensional submanifold.

The following 2-form equips $X$ with a near-symplectic structure.

$$\omega = dt\wedge df + \ast(dt\wedge df) =: \omega_{A} + \omega_{B}$$ 

Here $\ast$ denotes the Hodge $\ast$-operator, which is defined with respect to the product metric on $S^1 \times Y$.  The singular locus $Z_{\omega} = \left\{ p \in X \mid \omega_p=0 \right\}$ is in this case  $S^1\times \textnormal{Crit}(f)$.

An underlying near-symplectic structure $\omega$ on $X$ determines a decomposition of the normal bundle of the singular circles into two subbundles.  One of these subbundles is of rank 1.   A component of $Z_{\omega}$ over which this line bundle is topologically trivial is called {\it even}, and the one for which it is non-trivial is called {\it odd} \cite{P07}. These correspond to the cases when the normal bundle of the fold singularity is orientable or not (respectively).

In \cite{ADK05} a close relationship between near-symplectic $4$--manifolds and BLFs (under the name of singular Lefschetz fibrations) was established. When $X$ is compact and $b_{2}^{+}>0$, Theorem 1 in  \cite{ADK05} proves that there exists a BLF whose fold singularities $\Gamma$ correspond exactly to the singular locus of a near-symplectic form. So we may now endow such an $X$ with a  singular rank 2 Poisson structure, following the construction in the proof of our main theorem.

As remarked in the introduction---and computed explicitly in the previous section---the symplectic forms $\omega_{\Sigma_{p}}$ associated to the Poisson bi-vector $\Pi$ tend to $\infty$ as they approach the singular sets $C$ and $\Gamma$ of the BLF in question. In contrast, the near-symplectic form $\omega$ is positive on $C$ and vanishes on $\Gamma$.  

\subsection{Connected sums with $S^2\times S^2$ and complex projective planes}

Haya-no provides explicit examples of broken Lefschetz fibrations on 4-manifolds including connected sums with copies of $S^2\times S^2$, and complex projective spaces, among others \cite{Ha11}. These are examples of BLFs with fibres 2-spheres or 2-tori, which are known in the literature as genus-1 broken Lefschetz fibrations. 

In particular, the manifold $S^2\times S^2$ has a broken Lefschetz fibration $f\colon S^{2}\times S^{2} \rightarrow S^{2}$ with two Lefschetz type singularities and one circle of folds.  The connected sum $ S^{2}\times S^{2} \#  S^{2}\times S^{2}$ is the total space of a BLF with four Lefschetz singularities and one singular circle.   For a more general case see \cite{Ha11}. Our main theorem equips these manifolds with Poisson structures whose generic symplectic leaves are precisely these  2-spheres and 2-tori fibres.

\subsection{4-sphere}
Our construction shows that the 4-sphere is an example of a 4-manifold that admits a Poisson structure of generic rank 2, but no near-symplectic structure, as $H^{2}(S^{4}) \simeq 0$.

 In the work of Auroux, Donaldson and Katzarkov it was shown that there is a singular fibration via a BLF $f\colon S^4  \rightarrow S^2$ \cite{ADK05}.  There is one singular circle that gets mapped to the equator of $S^{2}$. The total space gets a decomposition into three pieces that are glued together.  Over the northern hemisphere $D^{2}_{N} \subset S^{2}$ the fibres are 2-spheres and $f^{-1}(D_{N}^{ 2}) \simeq S^{2} \times D^{2}:= X_{+}$.  Fibres over the southern hemisphere $D^{2}_{S}$ are 2-tori and the $f^{-1}(D_{S}) \simeq T^{2} \times D^{2}:= X_{-}$. Over the equatorial strip $E=S^{1} \times I$ lies the product of the singular circle $S^{1}$ with the standard cobordism from $T^{2}$ to $S^{2}$, which is diffeomorphic to a solid torus with a ball removed. That is $f^{-1}(E) \simeq S^{1}\times ((S^{1} \times D^{2}) \setminus B^{3})$.  Near the circle of folds, a Poisson structure of rank 2 on associated to this fibration on $S^{4}$ can be described as in equation \eqref{eq:pi-gamma}.
On the regular regions $X_{+}$ and $ X_{-}$ the Poisson bivector gets defined by the symplectic form of $S^{2}$ and $T^{2}$.

\end{document}